\magnification=1200
\overfullrule=0pt
\centerline {{\bf  Existence of zeros for operators taking their values 
in the dual of a Banach space}}\par
\bigskip
\bigskip
\centerline {BIAGIO RICCERI}\par
\bigskip
\bigskip
\bigskip
\bigskip
Throughout the sequel, $E$ denotes a reflexive real Banach space and $E^*$ its
 topological dual. We also assume that $E$ is locally uniformly convex. This
 means that for each $x\in E$, with $\|x\|=1$, and each $\epsilon>0$ there
 exists $\delta>0$ such that, for every $y\in E$ satisfying $\|y\|=1$ and
 $\|x-y\|\geq \epsilon$, one has $\|x+y\|\leq 2(1-\delta)$. Recall 
 that any reflexive Banach space admits an equivalent norm with which it is
 locally uniformly convex ([1], p. 289).
 For $r>0$, we
set $B_{r}=\{x\in E : \|x\|\leq r\}$. \par
 \smallskip
 Moreover, we
fix a topology $\tau$ on $E$, weaker than the strong topology and stronger than
the weak topology, such that $(E,\tau)$ is
 a Hausdorff locally convex
topological vector space with the property that the $\tau$-closed convex hull
of any $\tau$-compact subset of $E$ is still $\tau$-compact and the
relativization of $\tau$
to $B_1$ is metrizable by a complete metric. In practice, the most usual choice
of $\tau$ is either the strong topology or the weak topology provided $E$ is
also separable.\par
\smallskip
The aim of this short paper is to establish the following result and present
 some of its consequences:\par
\medskip
THEOREM 1. - {\it Let $X$ be a paracompact topological space and $A:X\to E^*$
a weakly continuous operator.
 Assume that there exist a number
 $r>0$, a continuous function $\alpha:X\to {\bf R}$ satsfying
$$|\alpha(x)|\leq r\|A(x)\|_{E^*}$$
for all $x\in X$, a closed set $C\subset X$,  and a $\tau$-continuous
 function
$g:C\to B_{r}$ satisfying $$A(x)(g(x))=\alpha(x)$$ for all $x\in C$, in such
a way that, for every $\tau$-continuous
function $\psi:X\to B_r$ satisfying $\psi_{|C}=g$, there exists $x_{0}\in X$
such that $$A(x_{0})(\psi(x_{0}))\neq\alpha(x_{0})\ .$$
Then, there exists $x^{*}\in X$ such that $A(x^{*})=0$.}\par
\medskip
For the reader's convenience, we recall
that
a multifunction $F:S\to 2^V$, between topological spaces, is said to be lower
semicontinuous at $s_{0}\in S$ if, for every open set $\Omega\subseteq V$
meeting $F(s_{0})$, there is a neighbourhood $U$ of $s_{0}$ such that
$F(s)\cap \Omega\neq \emptyset$ for all $s\in U$. $F$ is said to be lower
semicontinuous if it is so at each point of $S$.\par
\smallskip
The following well-known results will be our main tools.\par
\medskip
THEOREM A ([3]). - {\it Let $X$ be a paracompact topological space
and $F:X\to 2^{B_1}$ a $\tau$-lower semicontinuous multifunction
with nonempty $\tau$-closed convex values.\par
Then, for each closed set $C\subset X$ and each $\tau$-continuous
function $g:C\to B_1$ satisfying $g(x)\in F(x)$ for all $x\in C$, there
exists a $\tau$-continuous function $\psi:X\to B_{1}$ such that
$\psi_{|C}=g$ and $\psi(x)\in F(x)$ for all $\in X$}.
\medskip
THEOREM B ([4]). - {\it Let $X, Y$ be two topological spaces,
with
$Y$ connected and locally connected, and let $f:X\times Y\to {\bf R}$
be a function satisfying the following conditions:\par
\noindent
$(a)$\hskip 10pt for each $x\in X$, the function $f(x,\cdot)$ is continuous,
 changes sign in $Y$
  and is identically zero in no nonempty open subset
   of\hskip 5pt $Y$ \ ;\par
  \noindent
$(b)$\hskip 10pt the set $\{(y,z)\in Y\times Y : \{x\in X:
f(x,y)<0<f(x,z)\}$ is open in $X\}$ is dense in $Y\times Y$\ .\par
Then, the multifunction $x\to \{y\in Y : f(x,y)=0$ and $y$ is not a
local extremum for $f(x,\cdot)\}$ is lower semicontinuous, and its values
are nonempty and closed.}\par
\medskip
{\bf Proof of Theorem 1.}  Arguing by contradiction, assume that $A(x)\neq 0$
for all $x\in X$. For each $x\in X$, $y\in B_1$, put
$$f(x,y)=A(x)(y)-{{\alpha(x)}\over {r}}$$
and
$$F(x)=\{z\in B_{1} : f(x,z)=0\}\ .$$
Also, set
$$X_{0}=\{x\in X : |\alpha(x)|<r\|A(x)\|_{E^*}\}\ .$$
Since $A$ is weakly continuous, the function
$x\to \|A(x)\|_{E^*}$, as supremum of a family of
continuous functions, is lower semicontinuous.
 From this, it follows that the set
$X_0$ is open. For each $x\in X_0$, the function
$f(x,\cdot)$
is continuous and has no local, nonabsolute, extrema, being affine.
Moreover, it changes sign in $B_1$ since $A(x)(B_{1})=
[-\|A(x)\|_{E^*},\|A(x)\|_{E^*}]$ (recall that $E$ is reflexive). Since
$f(\cdot,y)$ is continuous for all $y\in B_1$, we then realize that the
restriction of $f$ to $X_{0}\times B_1$ satisfies the hypotheses of
Theorem B, $B_1$ being considered with the relativization of the
strong topology.
 Hence, the multifunction $F_{|X_{0}}$ is lower semicontinuous.
Consequently, since $X_{0}$ is open, the multifunction $F$ is lower
semicontinuous at each point of $X_{0}$. Now, fix $x_{0}\in X\setminus
X_{0}$. So, $|\alpha(x_{0})|=r\|A(x_{0})\|_{E^*}$.
Let $y_{0}\in F(x_{0})$ and $\epsilon>0$. Clearly, since $y_{0}$ is
an absolute extremum of $A(x_{0})$ in $B_1$, one has
$\|y_{0}\|=1$.
Choose $\delta>0$ so that,
for each $y\in E$ satisfying $\|y\|=1$ and $\|y-y_{0}\|\geq \epsilon$,
one has $\|y+y_{0}\|\leq 2(1-\delta)$.
By semicontinuity, the function
$x\to (\|A(x)\|_{E^*})^{-1}$ is bounded in some neighbourhood of $x_0$,
and so, since the functions $\alpha$ and $A(\cdot)(y_{0})$ are continuous,
it follows that
$$\lim_{x\to x_{0}}{{\left | A(x)(y_{0})-{{\alpha(x)}\over {r}}\right |}\over
{\|A(x)\|_{E^*}}}=0\ .$$
So, there is a neighbourhood $U$ of $x_0$ such that
$${{\left | A(x)(y_{0})-{{\alpha(x)}\over {r}}\right |}\over
{\|A(x)\|_{E^*}}}<{{\epsilon\delta}\over {2}} \eqno{(1)}$$
for all $x\in U$. Fix $x\in U$. Pick $z\in E$, with $\|z\|=1$, in such a
way that 
$|A(x)(z)|=\|A(x)\|_{E^*}$ and
$$\left ( A(x)(z) - {{\alpha(x)}\over {r}}\right )
\left ( A(x)(y_{0}) - {{\alpha(x)}\over {r}}\right ) \leq 0\ .$$
>From this choice, it follows, of course, that the segment joining $y_{0}$ and
$z$ meets the hyperplane $(A(x))^{-1}({{\alpha(x)}\over {r}})$. In
other words, there is $\lambda\in [0,1]$ such that
$$A(x)(\lambda z + (1-\lambda)y_{0})={{\alpha(x)}\over {r}}\ . \eqno{(2)}$$ 
So, if we put
$y=\lambda z + (1-\lambda)y_{0}$, we have $y\in F(x)$ and 
$$\|y-y_{0}\|=\lambda\|z-y_{0}\|\ . \eqno{(3)}$$
We claim that $\|y-y_{0}\|<\epsilon$. This follows at once from $(3)$ if
$\lambda<{{\epsilon}\over {2}}$. Thus,
 assume $\lambda\geq {{\epsilon}\over
{2}}$. In this case, to prove our claim, it is enough to show
that
$$2(1-\delta)<\|z+y_{0}\| \eqno{(4)}$$
since $(4)$ implies $\|z-y_{0}\|<\epsilon$. To this end, note that,
by $(2)$, one has
$${{\left | A(x)(y_{0})-{{\alpha(x)}\over {r}}\right |}\over
{\|A(x)\|_{E^*}}}={{\lambda |A(x)(z-y_{0})|}\over {\|A(x)\|_{E^*}}}\ ,$$
and so, from $(1)$, it follows that
$${{|A(x)(z-y_{0})|}\over {\|A(x)\|_{E^*}}}<\delta\ . \eqno{(5)}$$
Suppose $A(x)(z)=\|A(x)\|_{E^*}$. Then, from $(5)$, we get
$$1-\delta<{{A(x)(y_{0})}\over {\|A(x)\|_{E^*}}}\ . \eqno{(6)}$$
On the other hand, we also have
$$1+{{A(x)(y_{0})}\over {\|A(x)|_{E^*}}}=
{{A(x)(z+y_{0})}\over {\|A(x)\|_{E^*}}}\leq \|z+y_{0}\|\ . \eqno{(7)}$$
So, $(4)$ follows from $(6)$ and $(7)$. Now, suppose $A(x)(z)=
-\|A(x)\|_{E^*}$. Then, from $(5)$, we get
$$1-\delta<-{{A(x)(y_{0})}\over {\|A(x)\|_{E^*}}}\ . \eqno{(8)}$$
On the other hand, we have
$$1-{{A(x)(y_{0})}\over {\|A(x)|_{E^*}}}=-
{{A(x)(z+y_{0})}\over {\|A(x)\|_{E^*}}}\leq \|z+y_{0}\|\ . \eqno{(9)}$$
So, in the present case, $(4)$ is a consequence of $(8)$ and $(9)$.
 In such a manner, we have proved that $F$ is lower
semicontinuous at $x_{0}$. Hence, it remains proved that
$F$ is lower semicontinuous in $X$ with respect
to the strong topology, and so, {\it a fortiori},
with respect to $\tau$. Since $F$ is also with
nonempty $\tau$-closed convex values, and ${{g}\over {r}}$ is a
$\tau$-continuous selection of it over the closed set $C$, by Theorem 
A, $F$ admits a $\tau$-continuous selection $\omega$ in
$X$ such that $\omega_{|C}={{g}\over {r}}$. At this point,
if we put $\psi=r\omega$, it follows that $\psi$ is a $\tau$-continuous 
function, from $X$ into $B_r$, such that $\psi_{|C}=g$ and
$A(x)(\psi(x))=\alpha(x)$ for all $x\in X$, against the hypotheses.
This concludes the proof.\hfill $\bigtriangleup$
\medskip
We now indicate two reasonable ways of application of Theorem 1. The
first one is based on the Tychonoff fixed point theorem.\par
\medskip
THEOREM 2. - {\it Assume that $E$ is a separable Hilbert space, with inner product
$\langle\cdot,\cdot\rangle$. Let $r>0$ and let $A:B_{r}\to E$ be a continuous
operator from the weak to the strong topology. Assume that there
exist a weakly continuous function $\alpha:B_{r}\to {\bf R}$
satisfying
$|\alpha(x)|\leq r\|A(x)\|_{E^*}$ for all $x\in B_r$, and a
weakly continuous function $g:C\to B_{r}$ such
$$\langle A(x),g(x)\rangle =\alpha(x)\hskip 5pt and\hskip 5pt  g(x)\neq x$$ for all
$x\in C$, where
$$C=\{x\in B_{r} : \langle A(x),x\rangle =\alpha(x)\}\ .$$
Then, there exists $x^{*}\in B_r$ such that $A(x^{*})=0$.}\par
\smallskip
PROOF. Identifying $E$ with $E^*$, we apply Theorem 1 taking $X=B_r$,
with the relativization of
the weak topology of $E$, and taking as $\tau$ the weak topology of
$E$. Due to the kind of continuity we are assuming for $A$, the
 function $x\to \langle A(x),x\rangle$ turns out to be weakly continuous
(see the proof of Theorem 4),
 and so the set $C$ is weakly closed. Now, let $\psi:B_{r}\to B_r$ be any
weakly continuous function such that $\psi_{|C}=g$. By the Tychonoff 
fixed point theorem, there is $x_{0}\in B_r$ such that $\psi(x_{0})=
x_0$. Since $g$ ha no fixed points in $C$, it follows that $x_{0}\notin
C$, and so $$\langle A(x_{0}),\psi(x_{0})\rangle=
\langle A(x_{0}),x_{0}\rangle \neq \alpha(x_{0})\ .$$
Hence, all the assumptions
of Theorem 1 are satisfied, and the conclusion follows from it.
\hfill $\bigtriangleup$\par
\medskip
It is worth noticing the following consequence of Theorem 2.\par
\medskip
THEOREM 3. - {\it Let $E$ and $A$ be as in Theorem 2. Assume that
for each $x\in B_r$, with $\|A(x)\|>r$, one has
$$\left \|A\left ( {{rA(x)}\over {\|A(x)\|}}\right ) \right \| \leq r\ .
\eqno{(10)}$$
Then, the operator $A$ has either a zero or a fixed point.}\par
\smallskip
PROOF. Define the function $\alpha:B_{r}\to {\bf R}$ by
$$\alpha(x)=\cases {\|A(x)\|^{2} & if $\|A(x)\|\leq r$\cr & \cr
r\|A(x)\| & if $\|A(x)\|>r\ .$\cr}$$
Clearly, the function $\alpha$ is weakly continuous and satisfies
$|\alpha(x)|\leq r\|A(x)\|$ for all $x\in B_r$. Put
$$C=\{x\in B_{r} : \langle A(x),x\rangle =\alpha(x)\}\ .$$
Note that if $x\in C$ then $\|A(x)\|\leq r$. Indeed, otherwise, we would have
$\langle A(x),x\rangle =r\|A(x)\|$, and so, necessarily, $x={{rA(x)}\over
{\|A(x)\|}}$, against $(10)$. Hence, we have $\langle A(x),A(x)\rangle =
\alpha(x)$ for all $x\in C$. At this point, the conclusion follows at once
from Theorem 2, taking $g=A_{|C}$.\hfill $\bigtriangleup$\par
\medskip
REMARK 1. - It would be interesting to know whether Theorem 3 can
be improved assuming that $A$ is a continuous operator with relatively
compact range.\par
\medskip
The second application of Theorem 1 is based on connectedness arguments.
For other results of this type we refer to [5] (see also [2]).\par
\medskip
THEOREM 4. - {\it Let $X$ be a connected paracompact topological space
and $A:X\to E^*$ a weakly continuous and locally bounded operator.
Assume that there exist
$r>0$, a closed set $C\subset X$, a continuous function $g:C\to B_r$
and an upper semicontinuous function $\beta:X\to {\bf R}$, with
$|\beta(x)|\leq r\|A(x)\|_{E^*}$ for all $x\in X$, such that
$g(C)$ is disconnected, $$\beta(x)\leq A(x)(g(x))$$ for all $x\in C$ and
$$A(x)(y)<\beta(x)$$ for all $x\in X\setminus C$ and $y\in B_{r}\setminus
g(C)$.\par
Then, there exists $x^{*}\in C$ such that $A(x^{*})=0$.}\par
\smallskip
PROOF. First, note that the function $x\to A(x)(g(x))$ is continuous in $C$.
To see this, let $x_{1}\in C$ and let $\{x_\gamma\}_{\gamma\in D}$
be any net in $C$
converging
to $x_{1}$. By assumption, there are $M>0$ and a neighbourhood $U$ of
$x_1$ such that $\|A(x)\|_{E^*}\leq M$ for all $x\in U$. Let $\gamma_{0}\in
D$ be such that $x_{\gamma}\in U$ for all $\gamma\geq \gamma_{0}$. Thus,
for each $\gamma\geq \gamma_{0}$, one has
$$|A(x_{\gamma})(g(x_{\gamma}))-A(x_{1})(g(x_{1}))|\leq
M\|g(x_{\gamma})-g(x_{1})\| + |A(x_{\gamma})(g(x_{1}))-A(x_{1})(g(x_{1}))|$$
from which, of course, it follows that
$\lim_{\gamma}A(x_{\gamma})(g(x_{\gamma}))=A(x_{1})(g(x_{1}))$. Next,
observe that the multifunction $x\to [\beta(x),r\|A(x)\|_{E^*}]$ is lower
semicontinuous and that the function $x\to A(x)(g(x))$ is a continuous
selection of it in $C$. Hence, by Michael's theorem, there is a continuous
function $\alpha:X\to {\bf R}$ such that $\alpha(x)=A(x)(g(x))$ for all
$x\in C$ and $\beta(x)\leq \alpha(x)\leq r\|A(x)\|_{E^*}$ for all $x\in
X$. Now, let $\psi:X\to B_r$ be any continuous function such that $\psi_{|C}=
g$. Since $X$ is connected, $\psi(X)$ is connected too. But then,
since $g(C)$ is
disconnected and $g(C)\subset \psi(X)$, there exists $y_{0} \in \psi(X)\setminus
g(C)$. Let $x_{0}\in X\setminus C$ be such that $\psi(x_{0})=y_{0}$.
So, by hypothesis, we have
$$A(x_{0})(\psi(x_{0}))=A(x_{0})(y_{0})<\beta(x_{0})\leq \alpha(x_{0})\ .$$
Hence, taking as $\tau$ the strong topology of $E$,  all the assumptions
of Theorem 1 are satisfied, and the conclusion follows from it.
\hfill $\bigtriangleup$\par
\medskip
REMARK 2. - Observe that when $X$ is first-countable, the local boundedness
of $A$ follows automatically from its weak continuity. This follows from the
fact that, in a Banach space, any weakly convergent sequence is bounded.\par
\medskip
It is worth noticing the corollary of Theorem 4 which comes out taking
$X=B_r$, $\beta=0$ and $g$=identity:
\par
\medskip
THEOREM 5. - {\it Let $E$ be a Hilbert space, with inner product
$\langle\cdot,\cdot\rangle$. Let
 $r>0$ and let $A:B_{r}\to E$ be
a continuous operator from the strong to the weak topology. Assume
that the set $C=\{x\in B_{r} : \langle A(x),x\rangle\geq 0\}$ is
disconnected and that, for each $x, y\in B_{r}\setminus C$, one
has $\langle A(x),y\rangle<0$.\par
Then, there exists $x^{*}\in C$ such that $A(x^{*})=0$.}\par
\bigskip
\bigskip
\noindent
{\bf Acknowledgement}. The author wishes to thank prof. J. Saint Raymond
for a useful correspondence.\par
\bigskip
\bigskip
\centerline {{\bf References}}\par
\bigskip
\bigskip
\noindent
[1]\hskip 10pt R. DEVILLE, G. GODEFROY and V. ZIZLER, {\it Smoothness
and renorming in Banach spaces}, Longman Scientific $\&$ Technical, 1993.\par
\smallskip
\noindent
[2]\hskip 10pt A. J. B. LOPES-PINTO, {\it On a new result on the existence of
zeros due to Ricceri}, J. Convex Anal., {\bf 5} (1998), 57-62.\par
\smallskip
\noindent
[3]\hskip 10pt E. MICHAEL, {\it A selection theorem}, Proc. Amer. Math. Soc.,
{\bf 17} (1966), 1404-1406.\par
\smallskip
\noindent
[4]\hskip 10pt  B. RICCERI, {\it Applications de th\'eor\`emes de
semi-continuit\'e inf\'erieure}, C. R. Acad. Sci. Paris, S\'erie I,
{\bf 295} (1982), 75-78.\par
\smallskip
\noindent
[5]\hskip 10pt B. RICCERI, {\it Existence of zeros via disconnectedness},
J. Convex Anal., {\bf 2} (1995), 287-290.\par
\bigskip
\bigskip
\bigskip
\bigskip
\bigskip
Department of Mathematics\par
University of Catania\par
Viale A. Doria 6\par
95125 Catania, Italy\par
{\it e-mail}: ricceri@dmi.unict.it
\bye